\documentclass{article}

\usepackage{arxiv}
\usepackage{multirow}
\usepackage[utf8]{inputenc} 
\usepackage[T1]{fontenc}    
\usepackage{hyperref}       
\usepackage{url}            
\usepackage{booktabs}       
\usepackage{amsfonts}       
\usepackage{nicefrac}       
\usepackage{microtype}      
\usepackage{lipsum}
\usepackage{verbatim}
\usepackage{comment}
\usepackage{graphics}
\usepackage{graphicx}
\usepackage{subfigure}
\usepackage{amssymb}
\usepackage{float}

\usepackage{setspace}
\usepackage{xcolor}
\usepackage{graphicx}
\usepackage{verbatim}
\usepackage{comment}
\usepackage{graphics}
 \usepackage{subcaption}
\usepackage{amsmath}
\usepackage{array,amsthm,amssymb,amsmath,enumitem}
\usepackage{amssymb}
\usepackage{amsthm} 
\usepackage{amsmath}
\usepackage{graphicx}
\usepackage{subcaption}
\usepackage{comment}
\usepackage{mathtools,nccmath}
\usepackage{siunitx}
\usepackage{graphicx}
\usepackage{verbatim}
\usepackage{comment}
\usepackage[ruled,norelsize]{algorithm}
\usepackage{algpseudocode,algorithmicx}
\usepackage[normalem]{ulem}
\newtheorem{theorem}{Theorem}

\newtheorem{definition}{Definition}
\usepackage{algorithm}
\DeclareMathOperator{\tr}{tr}
\DeclareMathOperator{\diag}{diag}

\title{Enhanced Error Bounds For The Masked Projection Techniques via Cosine-Sine Decompositionn}

\author{
 Brij Nandan Tripathi \\
  Department of Electronics and Electrical Engineering\\
  Indian Institute of Technology Guwahati\\
  Assam, India \\
  \texttt{brij18@iitg.ac.in} \\
   \And
 Hanumant Singh Shekhawat \\
   Department of Electronics and Electrical Engineering\\
  Indian Institute of Technology Guwahati\\
  Assam, India \\
  \texttt{h.s.shekhawat@iitg.ac.in} \\ 
}

\begin{document}
\maketitle
\begin{abstract}
The masked projection techniques are popular in the area of non-linear model reduction. Quantifying and minimizing the error in model reduction, particularly from masked projections, is important. The exact error expressions are often infeasible. This leads to the use of error-bound expressions in the literature.
In this paper, we derive two generalized error bounds using cosine-sine decomposition for uniquely determined masked projection techniques. Generally, the masked projection technique is employed to efficiently approximate non-linear functions in the model reduction of dynamical systems. The discrete empirical interpolation method (DEIM) is also a masked projection technique; therefore, the proposed error bounds apply to DEIM projection errors. Furthermore, the proposed error bounds are shown tighter than those currently available in the literature.
\end{abstract}


\section{Introduction}
Consider the following linear parameterized non-linear dynamical system 
\begin{align}\label{sys}
\frac{d}{dt} x(t; \mu) = f(x(t; \mu)); \quad x(t_0; \mu) = x_0;
\end{align}
where \( x(t; \mu) \in \mathbb{R}^n \) represents the state vector of dimension $n$, $\mu=\begin{bmatrix}\mu_1,\mu_2,\dots, \mu_N \end{bmatrix}^T \in \mathbb{S} \subseteq \mathbb{R}^N$ denotes the system parameters, and $f: \mathbb{R}^n \to \mathbb{R}^n$ is a nonlinear function. 
In many practical applications, reducing the dimension of \eqref{sys} from $n$ to a smaller dimension
$p \ll n$ is necessary. This process is known as model reduction \cite{antoulas2020interpolatory}.

A popular technique of projecting the states $x(t;\mu)$ onto a smaller dimension space does not help much due to the large number of nonlinear components of $f$ \cite{MPE, DEIM}.  A commonly used method to circumvent this problem is to select a few components of $f$, say $[f_{j_1},f_{j_2}\dots, f_{j_m}]$, and then interpolate to get an approximation $\tilde{f}$ (denoted by tilde over $f$) of $f$. The selected components of $f$ can be obtained  by a simple  \emph{selection/mask  operator}  $P$ such that $[f_{j_1},f_{j_2}\dots, f_{j_m}]^T=P^T f$.  Mathematically, we seek the best approximation of  $f$  (called masked projection) that is determined only from the component obtained by the selection operator $P$  and contained in $p$-dimensional space $\mathcal{U}$.  Usually,  $\mathcal{U}$  is given as the column space of a matrix $U_1$, where $U_1$ is obtained by data-driven approaches (see \cite{Antoulas,benner2015survey}). 

There are many methods based on masked projection in the area of model order reduction \cite{MPE,DEIM,QDEIM,LDEIM,galbally}.  The Gauss-Newton with an approximate tensor (GNAT) \cite{GNAT}, missing point estimation \cite{MPE}, and \cite{galbally}, have used the concept of a gappy POD. In contrast, the discrete empirical interpolation method (DEIM) \cite{DEIM} is interpolation-based. The discrete version of the empirical interpolation method \cite{EIM} is DEIM, which produces interpolation points from a given set of basis to approximate the nonlinear function. 

Quantifying the error in model reduction is essential to ensure that the reduced model approximates the original system within acceptable limits. The masked projection also adds significant errors in the model reduction. Hence, it is imperative to reduce this error. However, an exact expression for the error introduced by the masking operation is not always feasible. Therefore, error-bound expressions are commonly reported in the literature for the masked projection methods  \cite{QDEIM,DEIM}.

In this paper, we have derived two error bounds for masked projection methods, which are typically used to approximate vector-valued functions (VVFs), by utilizing the Cosine-Sine (CS) decomposition \cite{CS,CSWebsite}. We especially focus on the uniquely determined masked projection where the number of selected components from the nonlinear function $f$ is the same as the dimension of the projection space $U$ \cite{AcceleratedMPE}. 
The proposed error bounds are numerically tighter than the state-of-the-art (SOA) bound provided by QDEIM \cite{DEIM,QDEIM, AcceleratedMPE}. 
The error bound we derive applies to any method using the uniquely determined masked projection, thereby encompassing the error bound for the state-of-the-art DEIM method within its scope. A similar approach is taken in Section 2 of \cite{AcceleratedMPE}, where a minor generalization of the standard DEIM error bound (see Lemma 3.2 of \cite{DEIM}) is presented. However, in our proposed work, we provide a more comprehensive and detailed analysis of this topic, along with the proposed enhancements in the bound.

The paper is organized as follows. The introduction of the masked projection method and its error analysis are discussed in section \ref{sec:sec2} and section \ref{sec:sec4}, respectively. Section \ref{sec:sec5} presents the numerical experiment results where we compare the accuracy of our proposed bound with the state-of-the-art methods. Finally, section \ref{sec:sec6} concludes the paper.

\subsection{Notation and convention} 
The operator $\diag: \mathbb{R}^m \rightarrow \mathbb{R}^{m \times m}$ transforms a vector into a square matrix by assigning the vector's values to the diagonal entries of the matrix, while all other entries remain zero. 

\section{Projection and Masked Projection}\label{sec:sec2}
In this section, we outlined masked projection methods in general, for which we have derived tighter error bounds in this paper. 

Let $y \in \mathbb{R}^n$ represent an arbitrary vector, and let $U_1 \in \mathbb{R}^{n \times m}$ be a column-orthogonal matrix. We can consider $y$ as a state solution or another vector of interest obtained from a specific simulation. We also assume that columns of $U_1$ form a basis for a $m$-dimensional subspace $\mathcal{U} \in \mathbb{R}^n$. 
The best approximation to $y$ contained in $\mathcal{U}$ is the orthogonal projection \[\hat{y}=U_1U_1^Ty.\]
\emph{We will use the hat notation for the orthogonal projection of any vector in this paper.}
Now,  we define the mask/selection operator.
\begin{definition}[Mask/Selection Operator]
 A mask/selection operator $P \in \mathbb{R}^{n \times m}$  contains $m$ columns of identity matrix $\mathbb{I}_n \in \mathbb{R}^{n \times n}$ as its columns. In particular, 
 \begin{equation}
     P=\begin{bmatrix}
         e_{\phi_1} & \dots e_{\phi_m}
     \end{bmatrix}
 \end{equation}
 where $e_{\phi_i}$ is ${\phi_i}^{th}$ column of the identity matrix $\mathbb{I}_n \in \mathbb{R}^{n \times n}$.
\end{definition}

Let \( J = \{j_1, \ldots, j_m\} \subseteq \{1, \ldots, n\} \) be a subset of indices with the corresponding selection matrix \( P = (e_{j_1}, \ldots, e_{j_m}) \in \mathbb{R}^{n \times m} \). We now focus on an approximation $\tilde{y}$ of $y$ given only a few components of $y$, say $\{y_{j_k}\}_{1\le k  \le m }$. In this context, $m$ typically refers to the number of interpolation points or selected points used for obtaining an approximation. Note that  $[y_{j_1},y_{j_2}\dots, y_{j_m}]^T=P^T y$.  
Furthermore,   the approximation of $y$ that is determined only from the  component obtained by  the selection operator $P$  and contained in $m$-dimensional space $\mathcal{U}$ is given by
\begin{align}\label{mpdef}
   \Tilde{y}=U_1(P^TU_1)^{-1}P^Ty.
\end{align}
Here,  $P^TU_1$ is assumed invertible \cite{DEIM}. \emph{The above $\tilde{y}$ is called the masked projection of $y$. We will use the tilde notation for the masked projection of any vector.}

Please note that we have not focused on the selection operator $P$ generation in this paper. If the mask operator $P$ is generated using the DEIM algorithm [4], the resulting masked projection is called the DEIM approximation.



\section{Error Analysis}\label{sec:sec4}
In this section, we derive the error associated with the masked projection of the vector-valued function. 

The DEIM \cite{DEIM} provides a method to construct the selection operator $P$ and illustrates the associated error bound. A tighter bound is presented in QDEIM method given in  \cite{QDEIM}. 
Furthermore, the error bound developed in the literature is specific to how the DEIM algorithm generates the selection operator $P$ and does not apply to the masked projection method in general. In this paper, we derive a tight error bound compared to available bounds in the literature. In addition, our bound applies to the general selection operator, meaning it is independent of the method used to generate the selection operator and applies to any masked projection method.


\begin{theorem}\label{thm:1}
Suppose columns of $U=\begin{bmatrix}
    U_1 & U_2
\end{bmatrix} \in \mathbb{R}^{n \times n}$ form a complete basis of  the range of the function $f:\mathbb{S} \rightarrow \mathbb{R}^n$ for all possible values of $\mu \in \mathbb{S} \subset{\mathbb{R}} $ where $U_1 \in \mathbb{R}^{n \times m}$, $U_1^TU_1=\mathbb{I}_m$, $U_2 \in \mathbb{R}^{n \times n-m}$ and $U_2^TU_2=\mathbb{I}_{n-m}$. Then, for the given selection operator $P \in \mathbb{R}^{n \times m}$ such that $(P^TU_1)^{-1}$ exists, the error between $f(\mu))$ and its  masked projection $\tilde{f}(\mu)$  satisfies
\begin{equation}\label{eqn6}
        \lVert f(\mu)-\Tilde{f}(\mu) \rVert_2^2 < \left(1 +  \sum_{i=1}^m \frac{1-\sigma_i^2}{\sigma_i^2}\right) \lVert(\mathbb{I}-U_1U_1^T)f(\mu) \rVert_2^2.
\end{equation}
where $\sigma_i$ is the $i^{th}$ singular value of $P^TU_1$. Moreover, the distance between the masked projection $\tilde{f}(\mu)$ and the orthogonal projection $\hat{f}=U_1U_1^Tf(\mu)$  of  $f(\mu)$ satisfies
\begin{align}\label{eqn234mm}
        \lVert \Tilde{f}(\mu)-\hat{f}(\mu) \rVert_2^2&=\sum_{i=1}^m \frac{1-\sigma_i^2}{\sigma_i^2} y_{i}^2 \\
        &\leq \sum_{i=1}^m \frac{1-\sigma_i^2}{\sigma_i^2} \lVert (\mathbb{I}-U_1U_1^T)f(\mu) \rVert_2^2 \notag
\end{align}
where $y:=V_2^TU_2^Tf$ with $V_2$ being the right singular matrix of $P^TU_2$. 

\end{theorem}

\begin{proof}
 Since columns of $U=\begin{bmatrix}
    U_1 & U_2
\end{bmatrix}$ form a complete basis of  the range of $f(\mu)$ for all possible values of $\mu$. Then, for any $\mu \in \mathbb{S} $,
\begin{equation}\label{eq412}
    f(\mu)=Uc(\mu)
\end{equation}

 For any masking or selection operator $P\in \mathbb{R}^{n \times m}$ such that $(P^TU_1)^{-1}$ exists, the approximation is as follows (section 2.2 of \cite{DEIM})
\begin{equation}
   \Tilde{f}(\mu)= U_1(P^TU_1)^{-1}P^Tf(\mu)
\end{equation}
Define ${f}_\mu:={f}(\mu)$ , $\Tilde{f}:=\Tilde{f}(\mu)$ and $\hat{f}:=\hat{f}(\mu)$. Now, using \eqref{eq412}, we get
\begin{equation}
    \begin{split}
       \Tilde{f}&=  U_1(P^TU_1)^{-1}P^T\begin{bmatrix}
           U_1 & U_2
       \end{bmatrix}\begin{bmatrix}
           k_1\\
           k_2
       \end{bmatrix} \\
             &= U_1(P^TU_1)^{-1}\begin{bmatrix}
           P^TU_1 & P^TU_2
       \end{bmatrix}\begin{bmatrix}
           k_1\\
           k_2
       \end{bmatrix} \\
            &=U_1\begin{bmatrix}
           \mathbb{I} & (P^TU_1)^{-1}P^TU_2
       \end{bmatrix}\begin{bmatrix}
           k_1\\
           k_2
       \end{bmatrix} \\
            &=U_1k_1+(U_1(P^TU_1)^{-1}P^T)U_2k_2
    \end{split}
\end{equation}
where $k_1:=U_1^Tf_\mu$ and $k_2:=U_2^Tf_\mu$. Now,
\begin{equation}
    \Tilde{f}=U_1U_1^Tf_\mu +U_1(P^TU_1)^{-1}P^TU_2U_2^Tf_\mu
\end{equation}
Further, let $\hat{f}$ denote the orthogonal projection of $f$ onto $U_1$. Then,
\begin{equation}\label{eq413}
    \Tilde{f}-\hat{f}=U_1(P^TU_1)^{-1}P^TU_2U_2^Tf_\mu
\end{equation}
Assume that $\begin{bmatrix}P & \breve{P} \end{bmatrix}$ is a permutation matrix. Now, we can construct an orthonormal matrix as,
\begin{align}
    Q&=\begin{bmatrix}P^T \\ \breve{P}^T \end{bmatrix}U=\begin{bmatrix}
        P^TU_1 & P^TU_2\\
       \breve{P}^TU_1 & \breve{P}^TU_2
    \end{bmatrix} \notag \\
    &=:\begin{bmatrix}
        Q_{11} & Q_{12}\\
        Q_{21} & Q_{22}
    \end{bmatrix}
\end{align}
Further, the cosine-sine decomposition of $Q$ can be written as \cite{CS,CSWebsite},
\begin{equation}\label{eqnm1}
\begin{bmatrix}
    Z_1^T & 0\\
    0     & Z_2^T
\end{bmatrix}\begin{bmatrix}
        Q_{11} & Q_{12}\\
        Q_{21} & Q_{22}
    \end{bmatrix}\begin{bmatrix}
    V_1 & 0\\
    0     & V_2
\end{bmatrix}=\begin{bmatrix}
    C & S &0\\
    -S & C & 0\\
    0 & 0 & \mathbb{I}_{n-2m}
\end{bmatrix},
\end{equation}
where $Z_1,Z_2,V_1$ and $V_2$ are orthogonal matrices, $C=\diag(\{\sigma_i\}_{i=1,\ldots,m})$ where $\sigma_i$ denotes the $i^{th}$ singular value of $P^TU_1$ 
and $S= \diag(\{\sqrt{1-\sigma_i^2}\}_{i=1,\ldots,m})$. Here, we assume that $\sigma_1\ge \sigma_2\ge \cdots\ge \sigma_m$. Note that $P^TU_1$ is assumed to be invertible, which implies $\sigma_i\ne0$ .
%
$Q_{11}=P^TU_1=Z_1CV_1^T$ and $Q_{12}=P^TU_2=Z_1\begin{bmatrix}
    S & 0
\end{bmatrix}V_2^T$. Hence, 
\begin{equation}\label{eq414}
    (P^TU_1)^{-1}(P^TU_2)=V_1C^{-1}\begin{bmatrix}
    S & 0
\end{bmatrix}V_2^T
\end{equation}
From \eqref{eq413} and \eqref{eq414},
\begin{align}\label{eqm1}
   \lVert \Tilde{f}-\hat{f} \rVert_2^2&=\lVert U_1(P^TU_1)^{-1}P^TU_2U_2^Tf_\mu \rVert_2^2 \notag\\
                                        &=(U_1V_1C^{-1}\begin{bmatrix}
    S & 0
\end{bmatrix}V_2^TU_2^Tf_\mu)^T \notag\\
&\quad \quad (U_1V_1C^{-1}\begin{bmatrix}
    S & 0
\end{bmatrix}V_2^TU_2^Tf_\mu) \notag\\
    &=f_\mu^TU_2V_2\begin{bmatrix}
    S \\
    0
\end{bmatrix}C^{-2}\begin{bmatrix}
    S & 0
\end{bmatrix}V_2^TU_2^Tf_\mu\\
    &=k_2^TV_2\begin{bmatrix}
    S \\
    0
\end{bmatrix}C^{-2}\begin{bmatrix}
    S & 0
\end{bmatrix}V_2^Tk_2 \notag
\end{align}
Define $\beta:=\lVert V_2^Tk_2 \rVert_{\infty}$. Now,
\begin{equation*}
    \begin{split}
        \lVert \Tilde{f}-\hat{f} \rVert_2^2 \leq \sum_{i=1}^m \frac{1-\sigma_i^2}{\sigma_i^2} \beta^2
    \end{split}
\end{equation*}
Since, $\beta \leq \lVert V_2^Tk_2 \rVert_{2} = \lVert k_2 \rVert_{2} =\lVert(\mathbb{I}-U_1U_1^T)f \rVert_2$, where the latter equality holds because $U=\begin{bmatrix}
    U_1 & U_2
\end{bmatrix}$ forms a complete basis for $f$, and the former equality $\lVert V_2^Tk_2 \rVert_{2}= \lVert k_2 \rVert_{2}$ holds because $V_2$ is an orthonormal matrix.
Using this, the bound can be given as

\begin{equation}\label{eqn789}
    \lVert \Tilde{f}-\hat{f} \rVert_2^2 \leq \sum_{i=1}^m (\sigma_i^{-2}-1) \lVert (\mathbb{I}-U_1U_1^T)f  \rVert_{2}^2
\end{equation}

Now, the bound on $\Tilde{f}$ with respect to $f_\mu$ can be given as
\begin{equation}\label{eqn7}
    \begin{split}
        \lVert f_\mu-\Tilde{f} \rVert_2^2 &= \lVert f_\mu-\Tilde{f}-\hat{f}+\hat{f} \rVert_2^2\\
                                      &=  \lVert f_\mu-\hat{f} \rVert_2^2 + \lVert \Tilde{f}-\hat{f} \rVert_2^2\\
                                      &\leq \lVert(\mathbb{I}-U_1U_1^T)f_\mu \rVert_2^2 \left(1 +  \sum_{i=1}^m (\sigma_i^{-2}-1)\right)\\
                                      &=: \mathcal{B}_{thm1P}
    \end{split}
\end{equation}
Functions $ f_\mu-\hat{f}=U_2U_2^Tf_\mu$ and $ \Tilde{f}-\hat{f}=U_1((P^TU_1)^{-1}P^Tf_\mu-U_1^Tf_\mu)$ are orthogonal as $f_\mu-\hat{f}$ is projection on $U_2$ while $\Tilde{f}-\hat{f}$ is projection on $U_1$.
\end{proof}

We have derived the bound for the masked projection $\tilde{f}(\mu)$ for a sample point $\mu$ in \ref{thm:1} (see \eqref{eqn6}). Define $f_i:=f(\mu_i)$. Consider $N$ sample points $\mu_1, \mu_2 \dots \mu_N$ from the domain $\mathbb{S}$ and corresponding samples $f_1, f_2, \dots f_N$ . 
Now,  extension of \eqref{eqn6} to $N$ samples is straightforward by taking an average:
\begin{align}\label{eqm2}
    \mathcal{E} &:= \frac{1}{N}\sum_{i=1}^N\lVert f(\mu_i)-\Tilde{f}(\mu_i) \rVert_2^2 \notag\\
      &\leq \frac{1}{N} \left(1 +  \sum_{j=1}^m \frac{1-\sigma_j^2}{\sigma_j^2}\right) \sum_{i=1}^N\lVert(\mathbb{I}-U_1U_1^T)f_i \rVert_2^2 
      =: \mathcal{B}_{thm1}
\end{align}
We have seen that \(\lVert \tilde{f}_i - \hat{f}_i \rVert\) is an essential component of \(\lVert f_\mu - \tilde{f} \rVert\) in the proof of Theorem \ref{thm:1}. Therefore, a tighter bound on \(\lVert \tilde{f}_i - \hat{f}_i \rVert\) (as obtained in the next theorem) will lead to a better bound on \(\lVert f_\mu - \tilde{f} \rVert\).

\begin{theorem}\label{thm2}
  Let $f:\mathbb{S} \to  \mathbb{R}^n$ be a nonlinear vector-valued function with $\mathbb{S} \in \mathbb{R}$.  Define $f_i:=f(\mu_i)$. Consider $N$ sample points $\mu_1, \mu_2 \dots \mu_N$ from the domain $\mathbb{S}$ and corresponding samples $f_1, f_2, \dots f_N$ then 
  \begin{align}
      \breve{\mathcal{E}}:=\frac{1}{N} \sum_{i=1}^N  \lVert \Tilde{f}_i- \hat{f}_i \rVert^2  &\leq  \frac{1}{N}\sum_{i=1}^m (\sigma_{m-i+1}^{-2}-1) \lambda_{i}(XX^T)
  \end{align}
 where $P$ and $U_i$ are defined in theorem \ref{thm:1}, $\sigma_i$ denote the $i^{th}$ singular value of $P^TU_1$, and $X:=U_2^T\begin{bmatrix}
  f_1 & f_2 & \dots & f_N  
\end{bmatrix}$. 
\end{theorem}

\begin{proof}
   Using \eqref{eqm1},  the square error of $\tilde{f}_i$ with respect to $\hat{f}_i$  can be written as below
\begin{equation}\label{eqn12341}
      \lVert \Tilde{f}_i- \hat{f}_i \rVert^2 = f_i^TU_2V_2\begin{bmatrix}
    S \\
    0
\end{bmatrix}C^{-2}\begin{bmatrix}
    S & 0
\end{bmatrix}V_2^TU_2^Tf_i
\end{equation}
where, $\Tilde{f}_i=U_1(P^TU_1)^{-1}P^Tf_i$ and $\hat{f}_i=U_1U_1^Tf_i$.
Now,  the  average square error can be written as
\begin{align}\label{eqn1m}
    \breve{\mathcal{E}}&=\frac{1}{N} \tr(  X^TV_2\begin{bmatrix}
    S \\
    0
\end{bmatrix}C^{-2}\begin{bmatrix}
    S & 0
\end{bmatrix}V_2^TX) \notag\\
     &=\frac{1}{N} \tr(\begin{bmatrix}
    S \\
    0
\end{bmatrix}C^{-2}\begin{bmatrix}
    S & 0
\end{bmatrix}V_2^TXX^TV_2) \notag\\
&\leq \frac{1}{N}\sum_{i=1}^m \lambda_i(S^2C^{-2}) \lambda_{i}(XX^T) 
\end{align}
Here, we have used \cite[Theorem B.1]{Marshall1979InequalitiesTO} in the last step. Note that
we have assumed $\lambda_1(A) \geq \cdots \geq \lambda_n(A) \geq 0$ for a given non-negative definite symmetric matrix $A$.

All singular values $\sigma_i$ of $P^TU_1$ are less than or equal to one \cite[Theorem 2.1]{Marshall1979InequalitiesTO}. Hence, the diagonal entries of $S^2C^{-2}$ ($\frac{1-\sigma_i^2}{\sigma_i^2}$)  follow $\frac{1-\sigma_i^2}{\sigma_i^2} \leq \frac{1-\sigma_{i+1}^2}{\sigma_{i+1}^2}$ as $\sigma_1\ge \sigma_2\ge \cdots\ge \sigma_m>0$. 
Using this fact,  the bound on the average error can be given as
\begin{equation}\label{eqn12342}
     \breve{\mathcal{E}}  \leq \frac{1}{N}\sum_{i=1}^m (\sigma_{m-i+1}^{-2}-1) \lambda_{i}(XX^T)
\end{equation} 
\end{proof}



Using \eqref{eqn12342}, a bound on the average error between the original function and the function approximated using the selection operator can be derived as

\begin{align}\label{eqn12343}
\mathcal{E} &= \frac{1}{N}\sum_{i=1}^N\lVert f(\mu_i)-\Tilde{f}(\mu_i) \rVert_2^2 \notag\\
 & \leq \frac{1}{N}\sum_{i=1}^N \lVert(\mathbb{I}-U_1U_1^T)f_i \rVert_2^2 + \frac{1}{N}\sum_{i=1}^m (\sigma_{m-i+1}^{-2}-1) \lambda_{i}(XX^T)\notag\\
 &=: \mathcal{B}_{thm2}
\end{align}





In the next section, we have shown numerically that the above bound $\mathcal{B}_{thm2}$ is tighter than $\mathcal{B}_{thm1}$ (see \eqref{eqm2}).

\subsection{Comparison with state-of-the-art  bound}

In this section, we first show that the proposed bound \eqref{eqn12342} is tighter than the state-of-the-art (SOTA) bound given in \cite{QDEIM}.
  
We first show that the proposed bound in \ref{thm2} is tighter than the SOTA bound. The SOTA bound is obtained using the following expression (refer to section 2.1 of \cite{QDEIM}):  
\begin{equation}\label{eq2mm}
    \lVert f - \Tilde{f} \rVert_2^2 \leq \lVert (P^T U_1)^{-1} \rVert_2^2 \lVert (\mathbb{I} - U_1 U_1^T)f \rVert_2^2
\end{equation}

The SOA bound is established by analyzing $\lVert (P^T U_1)^{-1} \rVert_2$, which, in our notation, corresponds to $\sigma_m^{-1}$, where $\sigma_m$ represents the $m^{\text{th}}$ (and the smallest) singular value of $P^T U_1$. We show that the bound in \eqref{eqn12342} is better than the bound in \eqref{eq2mm}, hence the SOTA bound.

Suppose $O=\begin{bmatrix}
    h_m & h_{m-1} \dots h_1 & h_n & h_{n-1} \dots h_{n-m}
\end{bmatrix}$ where $h_i$ denotes the $i^{th}$ eigenvector of $XX^T$ (with assumption $\lambda_{i+1}(XX^T)\geq\lambda_i(XX^T)$ where $\lambda_i(XX^T)$ is $i^{th}$ eigenvalue of $XX^T$). Then, $O^TXX^TO=D$ where $D$ is diagonal matrix with $d_{ii}=\lambda_{n-i}(XX^T)$. Note that $O$ is an orthognal matrix as it contains the eigenvector of symmetric matrix. The error $\breve{\mathcal{E}}$ with $V_2=O$ can be written as
\begin{align}
    \breve{\mathcal{E}}
     &=\frac{1}{N} \tr(\begin{bmatrix}
    S \\
    0
\end{bmatrix}C^{-2}\begin{bmatrix}
    S & 0
\end{bmatrix}O^TXX^TO) \notag\\
&=\frac{1}{N} \tr(\begin{bmatrix}
    S^2C^{-2} & 0 \\
    0 & 0
\end{bmatrix}D)\\
&= \frac{1}{N} tr(S^2C^{-2}D[1:m,1:m])\\
&= \frac{1}{N}\sum_{i=1}^m (\sigma_{m-i+1}^{-2}-1) \lambda_{i}(XX^T)
\end{align}
where $D[1:m,1:m] \in \mathbb{R}^{m \times m}$ is principle submatrix of $D$. This implies that the equality of the proposed bound \ref{thm2} occurs at $O$.

 Since $O$ is an orthonormal matrix, $\lVert O^T U_2^T f \rVert = \lVert U_2^T f \rVert = \lVert (\mathbb{I} - U_1 U_1^T)f \rVert$. Suppose $(O^T U_2^T f)^T =
\begin{bmatrix}
    \beta_1 & \beta_2 \dots \beta_{n-m}
\end{bmatrix}^T$
Then the SOTA bound can be written as 
\begin{align}\label{eq2m}
    \lVert f-\Tilde{f} \rVert_2^2 &\leq \sum_{i=1}^{n-m} \sigma_m^{-2} \beta_i^2\\
                                    &= \lVert(\mathbb{I}-U_1U_1^T)f \rVert_2^2 + \sum_{i=1}^{n-m} (\sigma_m^{-2}-1) \beta_i^2
\end{align} 

Although the bound presented is bound on average error in \ref{thm2} but for $V_2=O$ and $N=1$, \eqref{eqn12343} can be written as 

\begin{align}
   e &=\frac{1}{N} \tr(\beta^T\begin{bmatrix}
    S \\
    0
\end{bmatrix}C^{-2}\begin{bmatrix}
    S & 0
\end{bmatrix}\beta^T) \\
&=\sum_{i=1}^m (\sigma_{m-i+1}^{-2}-1) \beta_i^2
\end{align}
Now, the bound on $\Tilde{f}$ with respect to $f_\mu$ can be given as
\begin{equation}\label{eqm3m}
    \begin{split}
        \lVert f_\mu-\Tilde{f} \rVert_2^2 &= \lVert f_\mu-\Tilde{f}-\hat{f}+\hat{f} \rVert_2^2\\
                                      &=  \lVert f_\mu-\hat{f} \rVert_2^2 + \lVert \Tilde{f}-\hat{f} \rVert_2^2\\
                                      &\leq \lVert(\mathbb{I}-U_1U_1^T)f_\mu \rVert_2^2 + \sum_{i=1}^m  (\sigma_{m-i+1}^{-2}-1) \beta_i^2\\
    \end{split}
\end{equation}

Each term ($\sigma_m^{-2} \beta_i^2$) in the sum of \eqref{eq2m} is greater or equal to the term ($\sigma_i^{-2}\beta_i^2$) of \eqref{eqm3m} as $\sigma_i \leq \sigma_{i+1}$ where $\sigma_i$ denotes the singular value $i^{th}$ of $P^TU_1$. This shows that the proposed bound is tighter than the SOA bound \cite{QDEIM}.

Next, we explain why the proposed bound in \ref{thm:1} yields better results than the SOA bound in numerical experiments.

In \eqref{eqn7}, let replace $\sigma_i$ with the max value which is $\lVert (P^TU_1)^{-1} \rVert_2^2 =\sigma_m^{-2}$ , 
\begin{equation}\label{eqnnn7}
    \begin{split}
        \lVert f_\mu-\Tilde{f} \rVert_2^2 
                                      \leq \left(1 +   m(\sigma_i^{-2}-1)\right) \lVert(\mathbb{I}-U_1U_1^T)f_\mu \rVert_2^2
    \end{split}
\end{equation}
Note that \eqref{eqn7} is a much tighter bound than the one presented in \eqref{eqnnn7}. Now the difference  between the bound presented in \eqref{eqnnn7} and SOA bound is defined as
\begin{equation*}
\begin{split}
    \alpha &= \left((1+m(n-m))-(1 +   m(\sigma_i^{-2}-1))\right) \lVert(\mathbb{I}-U_1U_1^T)f_\mu \rVert_2^2 \\
      &= m((n-m+1) - \sigma_m^{-2})\lVert(\mathbb{I}-U_1U_1^T)f_\mu \rVert_2^2
\end{split}    
\end{equation*}
It is evident that if $\alpha > 0$ then bound in \eqref{eqnnn7} is tighter than the SOA bound. Furthermore, $\alpha >0$ holds when $\lVert (P^TU_1)^{-1} \rVert_2^2 =\sigma_m^{-2} < n-m+1$.
In this regard, although the bound on $\lVert (P^TU_1)^{-1} \rVert_2^2 \leq 1+m(n-m)$ is provided in \cite{GOREINOV19971}, the authors acknowledged that the provided bound is not sharp and stated that they have not found any matrix for which the inequality $\lVert (P^TU_1)^{-1} \rVert_2^2 =\sigma_m^{-2} < n$ is violated. The same is acknowledged in \cite{QDEIM}. If this is the case, the bound in \eqref{eqnnn7}, although much looser than the one in \eqref{eqn7}, provides better result than the SOA bound.

\section{Numerical Experiments}\label{sec:sec5}
In this section, we demonstrate the tightness of the proposed error bounds through examples.

\subsection{Example 1}\label{sec:eg}
Consider function f : $\Omega \times D\rightarrow \mathbb{R}$
 \begin{equation*}
     f(x,y,\mu) =\frac{y}{\sqrt{(x+y-\mu)^2+(2x-3\mu)^2+0.01^2}}
  \end{equation*} 
where $(x,y) \in \Omega= [0,2]^2 \subseteq \mathbb{R}^2$ and $\mu \in D= [0,2]\subseteq \mathbb{R}$. We will first compute the masked projection $\tilde f$ of $f$ using \eqref{mpdef}. Recall that we have to use a selection operator $P$ such that $P^TU_1$  is invertible. We have used the DEIM method to generate such a $P$ \cite{DEIM}.  The generation of $U_1$ is explained as follows. Let $(x_i,y_j)$ be uniform grid on $\Omega$ for $i=1,\dots,n_x$ and $j=1,\dots , n_y$. Define $\mathbf{f}:D \rightarrow \mathbb{R}^{n_x \times n_y}$
 \begin{equation*}
     \mathbf{f}(\mu) = [f(x_i,y_j,\mu)] \in \mathbb{R}^{n_x \times n_y}
  \end{equation*} 
  for $\mu \in D$ and $i=1,\dots,n_x$ and $j=1,\dots , n_y$. 
The $225$ snapshots constructed from uniformly selected parameters from $D$ are used for constructing the POD basis $U_1$.
  The sampled data for different values of $\mu$ is arranged as a third-order tensor:
  \begin{equation*}
     X(:,:,k) = [\mathbf{f}(x_i,y_j,\mu_k)] \in \mathbb{R}^{n_x \times n_y \times n_{\mu}}
  \end{equation*}
 The tensor $X$ is unfolded along the $x$ direction into a matrix of dimension $(n_x*n_y) \times n_{\mu}$ to generate the empirical basis through the POD method \cite{Kolda}. We select the top $m$ (which denotes the number of interpolation points or number of selected points) to generate $U_1$. 
 
 A different set of $n_{\mu}=400$ equally spaced parameters in $D$ are used to measure the accuracy; the average numerical error is used and defined as,
  \begin{equation}\label{eqn999}
      \bar{\mathcal{E}}(f,\tilde f)=\frac{1}{n_{\mu}}\sum_{i=1}^{n_{\mu}} \lVert Y_{\mu_i}-\Tilde{Y}_{\mu_i} \rVert_2^2
  \end{equation}
  where  $Y_{\mu_k}:=vec(f(x_i,y_j,\mu_k))$  and $\tilde{Y}_{\mu_k}:=vec(\tilde{f}(x_i,y_j,\mu_k))$.
 Table \ref{tab:table1} and \ref{tab:table2} show the average numerical error $\bar{\mathcal{E}}(f,\tilde f)$ defined in \eqref{eqn999} for $(n_x=30,n_y=30,n_{\mu}=400)$ and $(n_x=40,n_y=40,n_{\mu}=400)$ respectively. It can be seen that the proposed error bounds $\mathcal{B}_{thm1}$ (see \eqref{eqm2}) and $\mathcal{B}_{thm2}$ (see  \eqref{eqn12343}) are more closer to $\bar{\mathcal{E}}(f,\tilde f)$ when compared to the bound provided in \cite{QDEIM} paper (See Theorem \ref{thm:thmQDEIM}). In fact, $\mathcal{B}_{thm2}$ is better than $\mathcal{B}_{thm1}$. Table \ref{tab:table1} and \ref{tab:table2} also show how the error and error bound vary on increasing the spatial discretization points. Since the QDEIM error bound is directly proportional to the number of spatial discretization points (see Theorem \ref{thm:thmQDEIM}), the error bounds in table \ref{tab:table2} increase significantly from those in table \ref{tab:table1}. However, this dependence is not there directly for the proposed bounds.

\begin{table}[!ht]
\centering
\resizebox{\columnwidth}{!}{\begin{tabular}{|c|c|c|c|c|}
\hline
\textbf{Number of} & \textbf{Average Error} & \textbf{Proposed Bound} & \textbf{Proposed Bound} & \textbf{Average Error} \\ 
\textbf{Interpolation Points(m)} & $\bar{\mathcal{E}}(f,\tilde f)$ using \textbf{\eqref{eqn999}} & \(\mathcal{B}_{\text{thm1}}\) in \eqref{eqm2} & \(\mathcal{B}_{\text{thm2}}\) in \eqref{eqn12343} & \textbf{Bound (QDEIM)} \\ \hline
4  & 1.088  & 387.95  & 150.1305  & 1744 \\ \hline
6  & 0.2190  & 92.20  & 26.42  & 563.32 \\ \hline
8  & 0.311  & 40.26   & 11.33  & 205.05 \\ \hline
10  & 0.023 & 5.733   & 1.96  & 80.22 \\ \hline
12 & 0.005 & 2.809   & 0.6188  & 31.73 \\ \hline
\end{tabular}}
\caption{A comparison of error corresponding to spatial discretization $(n_x,n_y)=(30,30)$ at various interpolation points for  Example 1}
\label{tab:table1}
\end{table}

\begin{table}[!ht]
\centering
\resizebox{\columnwidth}{!}{\begin{tabular}{|c|c|c|c|c|}
\hline
\textbf{Number of} & \textbf{Average Error} & \textbf{Proposed Bound} & \textbf{Proposed Bound} & \textbf{Average Error} \\ 
\textbf{Interpolation Points(m)} & $\bar{\mathcal{E}}(f,\tilde f)$ using \textbf{\eqref{eqn999}} & \(\mathcal{B}_{\text{thm1}}\) in \eqref{eqm2} & \(\mathcal{B}_{\text{thm2}}\) in \eqref{eqn12343} & \textbf{Bound (QDEIM)} \\ \hline
4  & 1.890  & 1242  & 485.95  & 5494 \\ \hline
6  & 0.5295  & 301.1  & 83.68  & 1779 \\ \hline
8  & 0.3302  & 126.38   & 34.58  & 655.58 \\ \hline
10  & 0.0447 & 28.05   & 5.99  & 266.65 \\ \hline
12 & 0.0217 & 10.19   & 2.02  & 114.97 \\ \hline
14 & 0.0067 & 4.23   & 0.8166  & 50.05 \\ \hline
\end{tabular}}
\caption{A comparison of error corresponding to spatial discretization $(n_x,n_y)=(40,40)$ at various interpolation points for  Example 1}
\label{tab:table2}
\end{table}

\subsection{Example 2}
The function f : $\Omega \times D\rightarrow \mathbb{R}$
\begin{align}
     f(x,\mu) &=(1-x) \cos{(3\pi\mu(x+1))}e^{-(x+1)\mu}
\end{align}
where $x$ $\in \Omega= [-1,1]$ and $\mu \in D= [1,\pi]\subseteq \mathbb{R}$. Let $x_i$ be uniform sampling on $\Omega$ for $i=1,\dots,n_x$. Define $\mathbf{f}:D \rightarrow \mathbb{R}^{n_x}$
 \begin{equation*}
     \mathbf{f}(\mu) = [f(x_i,\mu)] \in \mathbb{R}^{n_x}
  \end{equation*} 
for $\mu \in D$ and $i=1,\dots,n_x$.
The $50$ snapshots constructed from uniformly selected parameters from $D$ are used for constructing the POD basis $U_1$ for two spatial discretizations $n_x=100$ and $n_x=200$. 

A different set of $n_{\mu}=101$ equally spaced parameters in $D$ are used to measure the accuracy; the average numerical error is used and defined as,
  \begin{equation}\label{eqn991}
      \bar{\mathcal{E}}(\mathbf{f},\tilde{\mathbf{f}})=\frac{1}{n_{\mu}}\sum_{i=1}^{n_{\mu}} \lVert \mathbf{f}_{\mu_i}-\Tilde{\mathbf{f}}_{\mu_i} \rVert_2^2
  \end{equation}
 Table \ref{tab:table3} and \ref{tab:table4} show the average numerical error $ \bar{\mathcal{E}}(\mathbf{f},\tilde{\mathbf{f}})$ defined in \eqref{eqn991} for $(n_x=100,n_{\mu}=101)$ and $(n_x=200,n_{\mu}=101)$ respectively. It can be seen that the proposed error bounds $\mathcal{B}_{thm1}$ (see \eqref{eqm2}) and $\mathcal{B}_{thm2}$ (see  \eqref{eqn12343}) are more closer to $\bar{\mathcal{E}}(f,\tilde f)$ when compared to the bound provided in \cite{QDEIM} paper (see Theorem \ref{thm:thmQDEIM}). Here also,  $\mathcal{B}_{thm2}$ is better than $\mathcal{B}_{thm1}$. Table \ref{tab:table3} and \ref{tab:table4} also show how the error and error bound vary on increasing the spatial discretization points. Similar to the example 
 in Section \ref{sec:eg}, 
  the QDEIM error bound increases significantly with an increase in the spatial discretization points.

\begin{table}[!ht]
\centering
\resizebox{\columnwidth}{!}{\begin{tabular}{|c|c|c|c|c|}
\hline
\textbf{Number of} & \textbf{Actual Average Error} & \textbf{Proposed Bound} & \textbf{Proposed Bound} & \textbf{Average Error} \\ 
\textbf{Interpolation Points(m)} & $\bar{\mathcal{E}}(f,\tilde f)$ using \textbf{\eqref{eqn999}} & \(\mathcal{B}_{\text{thm1}}\) in \eqref{eqm2} & \(\mathcal{B}_{\text{thm2}}\) in \eqref{eqn12343} & \textbf{Bound (QDEIM)} \\ \hline
4  & 10.356  & 205.21  & 63.50  & 2441 \\ \hline
6  & 5.2671  & 137.25   & 38.09  & 1438 \\ \hline
8  & 2.2861 & 65.74   & 17.36  & 678.27 \\ \hline
10 & 0.6515 & 23.45   & 5.42  & 261.55 \\ \hline
14 & 00339 & 1.4305   & 0.2775  & 18.69 \\ \hline
\end{tabular}}
\caption{A comparison of error corresponding to spatial discretization $n_x=100$ at various interpolation points for example 2 }
\label{tab:table3}
\end{table}

\begin{table}[!ht]
\centering
\resizebox{\columnwidth}{!}{\begin{tabular}{|c|c|c|c|c|}
\hline
\textbf{Number of} & \textbf{Actual Average Error} & \textbf{Proposed Bound} & \textbf{Proposed Bound} & \textbf{Average Error} \\ 
\textbf{Interpolation Points(m)} & $\bar{\mathcal{E}}(f,\tilde f)$ using \textbf{\eqref{eqn999}} & \(\mathcal{B}_{\text{thm1}}\) in \eqref{eqm2} & \(\mathcal{B}_{\text{thm2}}\) in \eqref{eqn12343} & \textbf{Bound (QDEIM)} \\ \hline
4  & 20.8280  & 852.386  & 251.673  & 9990 \\ \hline
6  & 11.7972  & 600.48   & 163.35  & 5955 \\ \hline
10  & 1.2935 & 94.10   & 20.49  & 1108 \\ \hline
14 & 0.0677 & 6.722   & 1.4558  & 81.18 \\ \hline
16 & 0.0134 & 1.0375   & 0.2270  & 12.58 \\ \hline
\end{tabular}}
\caption{A comparison of error corresponding to spatial discretization $n_x=200$ at various interpolation points in example 2}
\label{tab:table4}
\end{table}

\section{Conclusion}\label{sec:sec6}
\label{section:sec5}
In this paper, we have derived two error bounds for uniquely determined masked projection and have shown its tightness with other state-of-the-art bounds using numerical examples.
We have also seen in our numerical experiment that the bound provided by theorem \ref{thm2}  ($\mathcal{B}_{thm2}$) is tighter than the bound provided by theorem \ref{thm:1} ($\mathcal{B}_{thm1}$).

\section{Appendix}
In this section, we have described the state or art error bound (QDEIM Error bound \cite{QDEIM}) in average error form $\mathcal{E}=\frac{1}{N} \sum_{i=1}^N \lVert f_i-\tilde{f}_i \rVert_2^2$ as below
\begin{theorem}(from \cite{QDEIM})\label{thm:thmQDEIM}
     Let $U_1 \in \mathbb{R}^{n \times m}$ be orthonormal ($U_1^TU_1 = I_m$, $m < n$). then there exists a masking operator $P$ such that the DEIM projection error is bounded by
     \begin{align}
       \frac{1}{N} \sum_{i=1}^N \lVert f_i-\tilde{f}_i \rVert_2^2 &\leq \frac{1}{N} (1+m(n-m))\sum_{i=1}^N \lVert f_i-U_1U_1^Tf_i \rVert_2^2
     \end{align}
\end{theorem}
where $\tilde{f}=U_1(P^TU_1)^{-1}P^Tf_i$




\end{document}